\documentclass[a4paper]{article}

\usepackage{booktabs}
\usepackage[english]{babel}
\usepackage[utf8]{inputenc}
\usepackage{amssymb}
\usepackage[fleqn]{amsmath}
\usepackage{arydshln}
\usepackage{amsthm}
\usepackage{graphicx}
\usepackage{framed}
\usepackage[colorinlistoftodos]{todonotes}
\usepackage{subcaption}
\usepackage{multirow}
\usepackage{enumitem}
\usepackage{titling} 
\usepackage{mathtools}
\usepackage{titlesec}
\usepackage[top=1.25in, bottom=1.25in, left=1.25in, right=1.25in]{geometry}
\usepackage{breqn}
\usepackage[numbers]{natbib}
\usepackage{tikz}
\usepackage{pgfplots}
\usepackage{ulem}
\usepackage{caption}
\usetikzlibrary{shapes,snakes}
\usepackage{stackengine}
\usetikzlibrary{arrows,%
	petri,%
	topaths}%
\usepackage{algorithm}
\usepackage{algpseudocode}
\usepackage{tabstackengine}

\TABbinary
\newcounter{example}%[chapter]
\newcounter{drafts}

\newcounter{module}
\newcounter{model}
\newcommand \linedabstractkw[2]{% Default width = 0.9
  \renewcommand\maketitlehookd{%
    \mbox{}\medskip\par
    \centering
    \hrule\medskip
    \begin{minipage}{0.9\textwidth}
    \textbf{Abstract}\\ #1\\
    
    \textit{Keywords: }#2
    \end{minipage}\medskip\hrule\medskip
    }      
}

\newcommand \ERPWtemplate[3]{
\usepackage{fancyhdr}
\pagestyle{fancy}
\fancyhf{}
\fancyhead[RO]{#1}
\fancyhead[LO]{#2}
\fancyfoot[CO]{#3}
\fancyfoot[RO]{\thepage}

\renewcommand{\headrulewidth}{1pt}
\renewcommand{\footrulewidth}{1pt}
}

\stackMath

\theoremstyle{plain}
\newtheorem{theorem}{Theorem}[section]

\theoremstyle{plain}
\newtheorem{lemma}[theorem]{Lemma}

\theoremstyle{plain}
\newtheorem{proposition}[theorem]{Proposition}

\theoremstyle{remark}

\theoremstyle{plain}

\theoremstyle{definition}

\theoremstyle{plain}

\ERPWtemplate{P. Wissing}{Self-converse mixed graphs}{\leftmark}

\tikzset{group/.style = {shape=circle,draw,dotted,minimum size=1em}}
\tikzset{vertex/.style = {shape=circle,draw,minimum size=1em}}
\tikzset{arc/.style = {->,> = latex'}}
\tikzset{edge/.style = {-,> = latex'}}
\tikzset{negarc/.style = {->,> = latex',dashed}}
\tikzset{negedge/.style = {-,> = latex',dashed}}
\tikzset{tree/.style = {-,> = latex',line width=.7mm}}

\title{
Self-converse mixed graphs are extremely rare
}
\author{
Pepijn Wissing\thanks{Corresponding author: p.wissing@tilburguniversity.edu}\\ \small{Department of Econometrics and Operations Research, Tilburg University}
}

\begin{document}
\linedabstractkw{
A  mixed graph is cospectral to its converse, with respect to the usual adjacency matrices. 
Hence, it is easy to see that a mixed graph whose eigenvalues occur uniquely, up to isomorphism, must be isomorphic to its converse. 
It is therefore natural to ask whether or not this is a common phenomenon.
This note contains the theoretical evidence to confirm that the fraction of self-converse mixed graphs tends to zero. 
}{
Mixed graph, Digraph, Self-converse 
}
\maketitle

\section{Introduction}
With the rising interest in spectral characterization of mixed graphs and some of their generalizations came an interesting question, concerning the existence of a fairly obvious pairs of cospectral mixed graphs. 
At the heart of this issue is the fact that a mixed graph and its converse, obtained from the former by reversing all of the oriented edges, are typically encoded by matrices that are each other's conjugate transpose. 
In other words, two mixed graphs that may not be equivalent, are almost trivially cospectral. 
Thus, in order for a mixed graph to be determined by its spectrum in the traditional way \cite{vandam2003}, it must be isomorphic to its converse; such mixed graphs are said to be \textit{self-converse} \cite{bondy2008graph}.
% the conjugate transpose of the matrix used to encode such an object almost trivially has the same eigenvalues, while such a matrix and its conjugate transpose do not necessarily describe an equivalent mixed graph

This then raises the following question: \textit{how rare are self-converse mixed graphs?}
In \cite{wissing2019negative}, numerical evidence (see Table \ref{tab: fraction self converse}, below) suggesting that the fraction of self-converse mixed graphs converges to zero as the number of vertices $n$ goes to infinity was provided, although a formal proof to this claim has not appeared yet. 
Specifically, while the counting polynomials by \cite{harary1966enumeration,harary1955number} are quite easily evaluated, they are relatively unwieldy objects to work with, for arbitrary $n$. 
In this note, we will present a simple proof, to formally show the desired result.

\section{Main result}
We recall some terminology.
Let $\Gamma=(V,E)$ be a graph with vertex set $V=\{1,\ldots,n\}$ and edge set $E\subseteq \binom{V}{2}$. 
A \textit{mixed graph} $X$ is obtained from $\Gamma$ by orienting each edge in $A\subseteq E(\Gamma)$ in some direction; the collection of undirected edges is denoted $E(X)$.
$\Gamma$ is said to be the \textit{underlying graph} of $X$, and the \textit{symmetric subgraph} $G(X)$ of $X$ is obtained by removing $A(X)$ from $X$.

Two (mixed) graphs $X$ and $Y$ are said to be \textit{isomorphic} if there exists a bijection $f: V(X)\to V(Y)$ such that $uv\in A(X)$ if and only if $f(u)f(v)\in A(Y)$, and $\{u,v\}\in E(X)$ if and only if $\{f(u),f(v)\}\in E(Y)$.
In case $X$ is mapped onto itself, $f$ is called an \textit{automorphism}. 
The \textit{converse} $X^c$ of $X$ is obtained from $X$ by reversing the direction of every arc in $A(X)$, and $X$ is said to be \textit{self-converse} if $X^c$ is equal to $X$, up to isomorphism.    
% $D$ digraph, $D^*$ converse

% $G$ symmetric subgraph
Finally, we recall the Erdös-Renyi random graph $\Gamma(n,p)$, and its natural mixed analog $X(n,p)$.
$\Gamma(n,p)$ is the order-$n$ graph such that every edge occurs with probability $p$. That is, $\mathbb{P}(\{u,v\}\in E)=p$.
Accordingly, $X(n,p)$ is the order-$n$ mixed graph whose arcs $uv$ occur with probability $p$; if both arcs $uv$ and $vu$ occur, we say instead that the edge $\{u,v\}$ occurs. 
Lastly, $X(n,p)$ is said to be \textit{asymmetric} if it has no non-identity automorphism. 

The key argument used in the proof of the main result is the notion that almost all symmetric subgraphs of a random mixed graph $X(n,1/2)$ have no nontrivial automorphism. 
For completeness, a proof of this essentially well-known fact for the desired Edr{\"o}s-Renyi graph $\Gamma(n,p=1/4)$ is included, below. 
For sufficiently large $n$, the following lemma should be clear.
\begin{lemma}\label{lemma assumption}
Let $\Gamma=\Gamma(n,1/4)$ and $\epsilon>0$ be arbitrarily small. 
For $n$ sufficiently large, the vertices of $\Gamma$ have degree at least $\frac{n}{4}(1-\epsilon)$  and at most $\frac{n}{8}(1+\epsilon)$ common neighbors, with high probability.
\end{lemma}
% \begin{proof}
% Note that $d_v\sim Bin(n,1/4)$ for every vertex $v$ of $G$. Then:
% \[\mathbb{P}\left(d_v\leq\frac{n}{4}(1-\epsilon)\right) \approx \mathbb{P}\left(Z\leq \frac{\frac{n}{4}(1-\epsilon) - \frac{n}{4}}{\sqrt{3n/16}}\right) = \mathbb{P}\left(Z\leq -\epsilon\sqrt{n}/\sqrt{3}\right)\to0 \text{ as }n\to\infty,\]
% where $Z\sim N(0,1).$
% \end{proof}
Now, the following is an easy adaptation from \cite[Thm. 3.1]{nevsetril2011graph}. 
\begin{theorem}[\cite{nevsetril2011graph}]\label{thm: assymetric random graph}
The probability that $\Gamma(n,1/4)$ is asymmetric tends to $1$ as $n\to\infty$. 
\end{theorem}
\begin{proof}
Let $V=\{1,2,\ldots,n\}$ be the vertex set of $\Gamma=\Gamma(n,1/4)$ and let $f: V\mapsto V$ be an automorphism such that $f(x)=y$ for some vertices $x\not=y$. 
Let $M=\{v\in V : f(v)\not=v\}$ be the set of vertices that are moved by $f$. 
Moreover, let $V'=\binom{V}{2},$ and let $f':V'\mapsto V'$ be the permutation defined by $f'(\{u,v\})=\{f(u),f(v)\}.$ 

By Lemma \ref{lemma assumption}, for sufficiently large $n$, there exist at least $\lceil\frac{n}{4}(1-\epsilon)-\frac{n}{8}(1+\epsilon)\rceil = \lceil\frac{n}{8}(1-3\epsilon)\rceil$ vertices that are connected by an edge to $x$ but not to $y$. 
All of these vertices are moved by the automorphism $f$. 
Therefore, $|M|\geq cn$ for $c=(1-3\epsilon)/8$ with $\epsilon$ small. 
Thus the number of pairs of vertices that are moved by this automorphism is at least $\binom{cn}{2}-n\geq c'n^2$ for 
\[ c'=\frac{c^2n-c}{2n}\xrightarrow[n\to\infty]{}\frac{1-6\epsilon+9\epsilon^2}{128}>0\text{ for } \epsilon\not=\frac13\]  and sufficiently large $n$.
Therefore, the number of cycles of $f'$ is at most $k=\binom{n}{2}-c'n^2/2$. 

If $f$ is an automorphism of $\Gamma$, then the pairs in one cycle of $f'$ are either all edges or they are all non-edges of $\Gamma$. 
Hence, there are at most $2^k$ graphs such that $f$ is their automorphism. 

Combining the above, it follows that the probability that $\Gamma(n,1/4)$ has a non-identity automorphism is at most %\todo{numerically, the number $c'$ is too small (also for graph analog)}most 
\[\frac{n!\cdot 2^{\binom{n}{2}-c'n^2/2}}{2^{\binom{n}{2}}}\leq \frac{n^n}{2^{c'n^2/2}},\]
which tends to $0$ as $n\to\infty$. 
Indeed, note that 
\[\log\left(\frac{n^n}{2^{c'n^2/2}}\right) = n\log n - \frac12c'n^2\log2 \xrightarrow[n\to\infty]{} -\infty,\]
for all $c'>0$.
\end{proof}
The next result now follows naturally, by observing that any relabeling of the vertices that maps a mixed graph $X$ to $X^c$ simultaneously maps its symmetric subgraph onto itself.
Indeed, since the latter implies with high probability that said mapping is, in fact, the identity mapping, a contradiction follows.  
\begin{proposition}\label{prop: main result}
The probability that $X(n,1/2)$ is is self-converse tends to zero as $n\to\infty$.
\end{proposition}
%  Dear Edwin,
%  now, that I (hopefully) understand the statement, I think there is a simple proof.
\begin{proof}
 Let $n\to\infty$, and let $X$ be an order-$n$ mixed graph whose symmetric subgraph is $G=G(X)$ %, that is, $G$ is the graph on the same vertex set with vertices $x$ and $y$ adjacent if and only if $D$ has both arrows $xy$ and $yx$.
 If $X=X(n,1/2)$, then $G$ is the Erd{\"o}s-Renyi graph with edge probability $\frac14$.
 By Theorem \ref{thm: assymetric random graph}, $G$ has no nontrivial automorphism with probability tending to $1$.  
 Now, since any isomorphism from $X$ to $X^c$ is an automorphism of $G$, said isomorphism must be the identity map.
 However, with a probability tending to $1$, there is a pair $(x,y)\in V\times V$ such that  $X$ contains the arc $xy$ but not its converse arc $yx$.
 Therefore, the identity map is no isomorphism from $X$ to $X^c$ (with high probability), thus yielding a contradiction.
 \end{proof}
%  Or do I still confuse anything?
One should be somewhat mindful of what is being counted. 
Proposition \ref{prop: main result} implies that the fraction of self-converse \textit{labeled} mixed graphs tends to zero, whereas we are interested in its unlabeled counterpart, i.e., the fraction of all non-isomorphic mixed graphs. 
Note the significant distinction: any mixed graph with only the identity automorphism has $n!$ labeled versions, whereas (e.g.) the complete graph only has one.
In other words, the former is weighted much more heavily than the latter, by a probabilistic argument. % many more times than the latter.
Fortunately, this does not invalidate the approach.
In their extensive book, \citet{harary1973graphical} prove that almost all graphs of order $n$ can be labeled in $n!$ ways, and observe:
\begin{theorem}[\cite{harary1973graphical}]
Most labeled graphs have property "P" if and only most unlabeled graphs have property "P".
\end{theorem}
It should be clear that the argumentation would directly carry over to mixed graphs.
Hence, the desired result follows from Proposition \ref{prop: main result}.
% Being mindful of what is being counted by the above, we may directly infer the desired conclusion. 
\begin{proposition}

\label{thm: main result}
The fraction of order-$n$ self-converse mixed graphs tends to zero as $n\to\infty$. 
\end{proposition}

\section{Convergence rate}
To give some idea as to the rate at which the fraction of self-converse mixed graphs tends to zero, we include Table \ref{tab: fraction self converse} from \cite{wissing2019negative}, below. 
Here, $f(n)$ denotes said fraction of the non-isomorphic mixed graph of order $n$, obtained by evaluation of counting polynomials from \cite{harary1966enumeration, harary1955number}.

\begin{table}[h!]
\begin{center}
\begin{tabular}{ccccccccccc}
$n$  & 3 & 4 & 5& 6& 7 &8 \\
$f(n)$ & 6.25$\cdot 10^{-1}$&   3.21$\cdot 10^{-1}$&   7.36$\cdot 10^{-2}$&   9.87$\cdot 10^{-3}$&   6.16$\cdot 10^{-4}$&   2.20$\cdot 10^{-5}$\\  
\hline
$n$ &9 &10  & 11 & 12 & 13 & 14  \\
$f(n)$ &   3.89$\cdot 10^{-7}$&   3.79$\cdot 10^{-9}$ & 1.85$\cdot 10^{-11}$&   4.89$\cdot 10^{-14}$&   6.50$\cdot 10^{-17}$&   4.58$\cdot 10^{-20}$\\
\hline
$n$  & 15& 16& 17 &18 &19 &20\\
$f(n)$ &   1.63$\cdot 10^{-23}$&   3.06$\cdot 10^{-27}$&   2.90$\cdot 10^{-31}$&   1.43$\cdot 10^{-35}$&   3.59$\cdot 10^{-40}$&   4.64$\cdot 10^{-45}$
  \end{tabular}
\caption{\label{tab: fraction self converse} The fraction $f(n)$ of mixed graphs of order $n$ that is self-converse.}
\end{center}
\end{table}

\section*{Acknowledgements}
The author would like to express his thanks to Oleg Verbitsky, for the useful observation that formed the core idea of this note. 
% Thanks to Oleg Verbitsky for the inspiration to this proof.%, and thanks to Edwin van Dam for his various helpful remarks. 

\scriptsize
\bibliography{mybib}{}
\bibliographystyle{plainnat}
\normalsize

\end{document}